\documentclass{article}
\usepackage{heron2e}
\usepackage{times} % if available; if not - comment out,
\title{Crossed Product Structure of Quantum
Euclidean Groups}
       \author{Tomasz Brzezi\'nski\footnote{Supported by the EPSRC
grant GR/K02244.} }
       \address{DAMTP, University of Cambridge
         Cambridge CB3 9EW, U.K.}
\def\sw#1{{\sb{(#1)}}}
\def\sco#1{{\sp{(\bar #1)}}}

\def\tens{\mathop{\otimes}}

\def\<{{\langle}}
\def\>{{\rangle}}

\def\id{{\rm id}} % roman id map
\def\eps{\epsilon}

\def\q2{{q^{-2}}}

\def\cross{\mbox{$\times\!\rule{0.3pt}{1.1ex}\,$}}
\def\note#1{{}}
\def\eqn#1#2{\begin{equation}#2\label{#1}\end{equation}}

\def\Z{{\bf Z}}
\def\R{{\bf R}}
\def\C{{\bf C}}
\begin{document}
\maketitle
\begin{abstract}
 It is shown
 that quantum Euclidean groups $E_q(2)$, $E_\kappa(2)$ and $E_\kappa(3)$ have 
 the structure of generalised crossed products.
\end{abstract}

\section{Introduction}
The idea of a crossed product, which appeared in studies of cohomology
of algebras over a Hopf algebra 
\cite{Swe:coh} and the Hopf-Galois theory  \cite{BlaCoh:cro}, 
 may be summarised as follows.
Given a Hopf algebra $H$ and an
algebra $M$ one wants to construct an algebra $P$ isomorphic to
$M\otimes H$ as a vector space and such that $P$ is a right
$H$-comodule with a coaction $\Delta_R :P\to P\otimes H$, which is an
algebra map trivial on $M$ and coinciding with
 a coproduct on $H$. The complete answer to this problem is known
 \cite{BlaCoh:cro}. 
 Namely, one can
construct such an algebra $P$ if and only if there is a weak action
$\rho : H\otimes M \to M$ of $H$ on $M$ and a {\it 2-cocycle} $\sigma
:H\tens H\to M$ which satisfy some conditions (see
\cite{BlaCoh:cro} for details).  

Crossed-products appear in a variety of  places in quantum group theory.
Hopf-Galois extensions  
are understood as quantum group principal bundles
\cite{BrzMa:gau} 
and a certain kind of a crossed product, known as a cleft extension,
corresponds to a trivial quantum group principal bundle. The quantum 
double of Drinfeld  is a special case of a crossed product, known as a double
cross product \cite{Maj:phy} (it is also an example of a trivial quantum
principal bundle). Also, $\kappa$-Poincar\'e
algebra \cite{LukNow:def} can be understood as a crossed product
\cite{MajRue:bic}. It seems reasonable  to expect that inhomogeneous
quantum groups which 
are built on a tensor product of two algebras should
enjoy the structure of some kind of a crossed product. One quickly
finds, however, that the standard definition of a crossed product is
too restrictive for this purpose. 
Take $E_q(2)$ \cite{VakKor:alg} as an
example. As a vector space, $E_q(2)$ is built on a tensor
product of the quantum hyperboloid $X_q$ \cite{Sch:int} and the Hopf algebra of
formal power series in $Z$, $\C [Z,Z^{-1}]$. Furthermore $E_q(2)$ is
a right $\C [Z,Z^{-1}]$-comodule  with a coaction that acts trivially
on $X_q$ and as a coproduct on $\C [Z,Z^{-1}]$. This coaction,
however, {\bf is not} an algebra map. Therefore, to interpret
$E_q(2)$, and some other inhomogeneous groups as well, 
as a crossed product, one needs to generalise the notion of a crossed
product in such a way that the condition that $\Delta_R$ be an algebra
map gets relaxed. Since this is the only place where the algebra
structure of $H$ enters, to achieve a suitable generalisation of the
crossed product we need to consider a crossed product by a {\bf
coalgebra} equipped with some additional structure, weaker then that
of a
Hopf algebra. This leads naturally to the concept of an
entwining structure which was introduced in \cite{BrzMa:coa} in an
attempt to construct a suitable quantum group gauge theory on quantum
homogeneous spaces. In this context a generalisation of a crossed
product  was proposed in \cite{Brz:cro}.

In this short note we summarise the construction of crossed products by a
coalgebra  of \cite{Brz:cro} and show that $E_q(2)$, $E_\kappa(2)$ and
$E_\kappa(3)$
are  such crossed products.

Throughout the paper we use standard Hopf algebra notations.
In a coalgebra $C$, $\Delta$ denotes the coproduct and
$\eps :C\to \C$ denotes the counit. We use the Sweedler notation to
denote the coproduct in $C$, $\Delta c = c\sw 1\tens c\sw 2$ (summation
understood), for any $c\in C$. By convolution product we mean a
product $*$ in a space of linear maps $C\to P$, where $P$ is an algebra,
given by $f*g(c) = f(c\sw 1)g(c\sw 2)$. A map $C\to P$ is said to be
convolution invertible if it is invertible with respect to this product.

\section{Crossed products by a coalgebra}
We first recall the definition of an entwining structure from
\cite{BrzMa:coa}.
We say that  a coalgebra $C$  and an algebra $P$ are {\em
entwined}
if there is a map  $\psi :C\otimes P\to
P\otimes C$ such that
\eqn{ent.A}{\psi\circ(\id_C\tens
\mu)=(\mu\tens\id_C)\circ\psi_{23}\circ\psi_{12},\quad \psi(c\tens
1)=1\tens   c,\quad \forall c\in C}
\[{(\id_P\tens\Delta)\circ\psi=\psi_{12}\circ\psi_{23}
\circ(\Delta\tens\id_P),\quad (\id_P\tens\eps)\circ\psi=\eps\tens\id_P,}\]
where $\mu$ denotes multiplication in $P$, and
$\psi_{12} = \psi\tens\id_P$ and $\psi_{23}=\id_P\tens\psi$. 
We denote the action of $\psi$ on $c\otimes u\in C\otimes P$ by
$\psi(c\otimes u)= u_\alpha\tens c^\alpha$ (summation understood).

Furthermore we assume that there is a group-like $e\in C$,
i.e. $\Delta e=e\tens e$, $\eps(e) =1$ and a map $\psi^C:C\tens C\to
C\tens C$ such that for any $c\in C$
$$
(\id\tens\Delta)\circ\psi\sp C=\psi\sp C_{12}\circ\psi^C_{23}
\circ(\Delta\tens\id), \quad
(\id\tens\eps)\circ\psi^C=\eps\tens\id,
\quad \psi^C(e\otimes c) = \Delta c,
$$
where $\psi_{12}^C = \psi^C\tens\id_C$ and
$\psi^C_{23}=\id_C\tens\psi^C$. We denote the 
action of $\psi^C$ on $b\tens c$ by $\psi^C(b\tens c) = c_A\tens b^A$
(summation understood).

With these assumptions $P$ is a right $C$-comodule with a coaction
$\Delta_R u = \psi(e\tens u)$. Moreover the fixed point subspace
$M=P^{coC}_e =\{ x\in P | \Delta_R x = x\tens e\}$ is a subalgebra of
$P$. We call $(P, C,\psi, e,\psi^C)$ the {\em entwining data}. A
number of examples of entwining data may be found in
\cite{BrzMa:coa}; all that are used in Section~3 come from the
structure of quantum homogenous spaces \cite{Brz:hom}. 

Given the entwining data $(P, C,\psi, e,\psi^C)$ and  
$M=P^{coC}_e$ the crossed product of $M$ by $C$ is defined as
follows. Assume that there are linear maps $\sigma :C\tens C\to
M$ and $\rho:C\tens P\to P$ such that for all $x,y\in M$,  $a,b,c\in C$:\\
(i) $\rho(e,x) = x, \quad \rho(c,1) = \eps(c)$;\\
(ii) $\rho(c\sw 1, (xy)_\alpha)\otimes c\sw 2 ^\alpha = \rho (c\sw
1,x_\alpha)\rho({c\sw 2}^\alpha\sw 1 , y_\beta)\otimes {c\sw
2}^\alpha\sw 2^\beta \in M\tens C$;\\
(iii) $\sigma(e, c) = \eps(c), \quad \sigma(c\sw 1, e_A)\tens {c\sw
2}^A = 1\tens c$;\\
(iv) $\rho(a\sw 1,\! \sigma(b\sw 1,\! c_A)_\alpha)\sigma(a\sw 2\!^{\alpha}\!\sw
1,\! b\sw 2\!^{A}_{~B})\tens a\sw 2\!^{\alpha}\!\sw 2\!^{B}\! = \!
\sigma(a\sw 1,\! b_A)\sigma(a\sw 2\!^{A}\!\sw 1,\! c_B)\!\tens a\sw
2\!^{A}\!\sw 2\!^{B}$;
(v) $\rho(a\sw 1,\! \rho(b\sw 1,\!x_\alpha)_\beta)\sigma(a\sw 2\!^{\beta}\!\sw
1,\! {b\sw 2\!^{\alpha}}_{A})\tens\! a\sw 2\!^{\beta}\!\sw
2\!^{A}\!\! = \!\!
\sigma(a\sw 1,\! b_A)\rho(a\sw 2\!^{A}\!\sw 1,\! x_\alpha)\!\tens\! a\sw
2\!^{A}\!\sw 2\!^{\alpha}$.

Then we can define an algebra structure on  $M\tens C$  with product
$$
(x\tens b)(y\tens c) = x\rho(b\sw 1,
y_\alpha)\sigma(b\sw{2}^{\alpha}\sw 1, c_A)\tens
b\sw{2}^{\alpha}\sw 2^{A}.
$$
and unit $1\tens e$.
This algebra is denoted by
$M\cross_{\rho,\sigma}C$ and called 
a {\em crossed product} by a coalgebra $C$. The pair $(\rho,\sigma)$ is called
the {\em crossed product data} for the entwining data $(P,C,\psi, e,
\psi^C)$.

Therefore, to construct a crossed product $M\cross_{\rho,\sigma}C$ one
first needs  to introduce  an entwining structure and then  maps
$\rho,\sigma$ that satisfy conditions (i-v). 
 One can easily
convince oneself  that all the previously known crossed products fit into above
scheme. For example, if one takes $C=H$ to be a Hopf algebra, $e=1$,
and $P$ a right $H$-comodule algebra with a coaction $\Delta_R :P\to
P\tens H$ 
then one can define a natural entwining structure by
$
\psi(h\tens u) = u\sco 0\tens hu\sco 1, \quad \psi^C(h\tens g) = g\sw
1\tens hg\sw 2,
$
where $\Delta_R u = u\sco 0\tens
u\sco 1$. With this entwining structure the conditions listed above
become the standard conditions for a weak action $\rho$ and a cocycle
$\sigma$ \cite{BlaCoh:cro}. 

Similarly, if 
     $C= B$ is a braided group, $P$ 
a right braided $B$-comodule algebra with   braiding $\Psi$, then we
take $e=1$ and define the entwining structure by 
$
\psi(b\tens u) = \Psi(b\tens u\sco 0)u\sco 1,
$
$
\psi^C(b\tens c)
= \Psi(b \tens c\sw 1)c\sw 2.
$ In this case the notion of a crossed product  by $B$ introduced
above, is  equivalent to the notion of a crossed
product by a braided group \cite{Ma:cro}.

A truly new crossed product is constructed assuming that 
 there is a convolution
invertible map 
$\Phi : C\to P$ such that,  for given 
entwining data $(P, C,\psi, e,\psi^C)$, $\Phi(e) =1$ and 
\begin{equation}
\psi\circ(\id_C\otimes\Phi) = (\Phi\otimes\id_C)\circ\psi^C.
\label{triv}
\end{equation}
In this case one can define the crossed product data $(\rho,\sigma)$ by
\begin{equation}
\rho(c,u) = \Phi(c\sw 1)u_\alpha\Phi^{-1}(c\sw 2^\alpha), \quad
\sigma(b,c) = \Phi(b\sw 1)\Phi(c_A)\Phi^{-1}(b\sw 2^A).
\label{cleft}
\end{equation}
This crossed product is isomorphic to $P$ as
an algebra. From the point of view of a coalgebra gauge theory
$P= M\cross_{\rho,\sigma}C$ is a trivial coalgebra bundle on $M$
\cite{BrzMa:coa}, while from the point of view of coalgebra
extensions it is a {\it cleft extension} of $M$. In the next section it will
be shown that $E_q(2)$, $E_\kappa(2)$ and
$E_\kappa(3)$ are crossed products of this type.

Given the entwining data $(P, C,\psi, e,\psi^C)$, crossed product data
$(\rho ,\sigma)$, and a convolution invertible map $\gamma : C\to M$, such that
$\gamma(e) =1$ and
$$ 
\psi^C_{23}\circ\psi_{12}\circ
(\id_C\tens\gamma\tens\id_C)\circ(\id_C\tens\Delta)
= (\gamma\tens\id_C\tens\id_C)\circ(\Delta\tens\id_C)\circ\psi^C,
$$
one can construct a new crossed product with data $(\rho',\sigma')$,
$$
\rho'(c,u) =\gamma(c\sw 1)\rho(c\sw 2, u_\alpha)\gamma^{-1}(c\sw 3^\alpha),
$$
$$
\sigma'(b,c) = 
\gamma(b\sw 1)\rho(b\sw 2, \gamma(c_A\sw 1)_\alpha)\sigma
(b\sw 3^\alpha, c_A\sw 2)\gamma^{-1}(b\sw 4^A).
$$
The crossed product data $(\rho,\sigma)$ and $(\rho',\sigma')$ lead to
the algebras which are isomorphic to each other and therefore
$(\rho,\sigma)$ and $(\rho',\sigma')$ are said to be  {\it gauge
equivalent}. From the point of view of the coalgebra gauge theory, the
map $\gamma$ is a gauge transformation of a trivial
principal bundle. The equivalence classes of $(\rho,\sigma)$ define a
``cohomology'' which generalises the non-Abelian cohomology introduced
in \cite{Ma:bey}.

\section{Crossed product structure of quantum Euclidean groups}

\subsection{The $q$-Euclidean group $E_q(2)$}
This example is
discussed in detail in \cite{Brz:cro}. The entwining data are as
follows. The coalgebra 
$C$ is spanned by group-like  $c_p$, $p\in {\bf Z}$. It can be
therefore equipped with an algebra structure of $\C [Z,Z^{-1}]$ with
$Z^p = c_{p+s}$ for any fixed $s\in \Z$. The total algebra $P$ is
taken to be  
$E_q(2)$ and is generated by $v,v^{-1}, n_+,n_-$ subject to the following 
relations
\begin{equation}
vn_\pm = q^2n_\pm v, \quad 
\quad n_+n_-
=q^2 n_-n_+ \quad vv^{-1} = v^{-1} v =1.
\label{eq2}
\end{equation}
As is well-known $E_q(2)$ is a quantum group, and thus it has a Hopf
algebra structure. For our purposes however, the Hopf algebra structure
of $E_q(2)$ is not important. We need the entwining structure on
$E_q(2)$ and $C$ and this is provided by the map 
$\psi:C\tens P\to P\tens C$  given by
$$
\psi(c_p\tens v^{\pm 1}) = v^{\pm 1}\tens c_{p\pm 1}, 
\quad \!\!\! \psi(c_p\tens n_\pm) = n_\pm\tens
c_p + \mu_\pm q^{2p}v^{\pm 1}\tens c_p -\mu_\pm q^{2p} v^{\pm 1}\tens
c_{p\pm 1} 
$$
where $\mu_+, \mu_-$ are non-zero complex numbers, and extended to the
whole of $E_q(2)$ by (\ref{ent.A}). It is an easy exercise to verify
that such an extension is compatible with (\ref{eq2}). Furthermore we
take $
e = c_s,$ and 
$
\psi^C(c_p\tens c_r) = c_r\tens c_{p+r-s}.
$

Given this entwining structure we can define a coaction $\Delta_R$ of $C$ on
$E_q(2)$, by $\Delta_R(u) = \psi(c_s\tens u)$. 
Following \cite{BonCic:fre} one easily finds that the fixed point
subalgebra $M$  of $P$ is 
generated by 
$
z_\pm = v^{\pm 1} + \mu_\pm^{-1}q^{-2s} n_\pm$,
which satisfy the relation
$
z_+z_- = q^2 z_- z_+ +(1-q^2).
$
Therefore $M$ is isomorphic to the quantum hyperboloid $X_q$ \cite{Sch:int}. 

To reveal the crossed
product structure of $E_q(2)$ we define a linear map 
$\Phi : C\to E_q(2)$ by
$\Phi(c_p) = v^{p-s}$. This map is clearly convolution invertible and
it satisfies (\ref{triv}). Therefore  $E_q(2)= X_q\cross_{\rho ,\sigma} C$
with a trivial cocycle $\sigma (c_p, c_r) =1$ and the map $\rho
:C\tens E_q(2)\to E_q(2)$,  
$$
\rho(c_p,  v^{\pm 1}) =1, \qquad 
\rho(c_p, n_\pm) = q^{2p}( q^{-2s}n_\pm +
\mu_\pm (v^{\pm 1}-1)).
$$

One can easily prove that, even if one introduces an algebra
structure on $C$ compatible with $\Delta$,  $E_q(2)$ is never a $C$-comodule
algebra. Furthermore, because $\psi(C\tens X_q)$ is {\it not} a subset
of $X_q\tens C$ there 
is no braiding making $E_q(2)$ a braided $C$-comodule
algebra. Therefore the notion of a coalgebra crossed product developed
in \cite{Brz:cro} and summarised  in this
paper is truly needed for description of internal structure of $E_q(2)$.

\subsection{The $\kappa$-Euclidean group $E_\kappa(2)$}
Similarly as for the $E_q(2)$ case, the coalgebra $C$
is spanned by group-like  elements $c_p$, $p\in {\bf Z}$. The
$\kappa$-deformation $E_\kappa(2)$ of two-dimensional Euclidean group,
which is obtained by the contraction of $SU_q(2)$ \cite{CelGia:...}
or, equivalently, by quantisation of the Poisson structure on $E(2)$
\cite{Mas:...}, is generated by $w,w^{-1}, a_1, a_2$, and its algebra
structure is determined by the relations 
$$
[w,a_1] ={\kappa\over 2}(w-1)^2, \quad [w,a_2]= \imath{\kappa\over
2}(w^2-1), \quad
[a_1,a_2] = \imath\kappa a_1.
$$
It is
clear that $a_1, a_2$ generate a subalgebra of $E_\kappa(2)$. This
subalgebra, known as a  {\it $\kappa$-plane}, is a homogeneous space of
$E_\kappa(2)$ \cite{BonCic:new}. We will denote it by $\R_\kappa^2$. 

We choose $e=c_0$ and proceed to define the entwining maps $\psi :
C\tens E_\kappa(2) \to 
E_\kappa(2)\tens C$ and $\psi^C:C\tens C\to C\tens C$. The former is
given by 
$$
\psi(c_p\tens w^{\pm 1}) = w^{\pm 1}\tens c_{p\pm 1}, \quad
%$$
%$$
\psi(c_p\tens a_\pm) = a_\pm\tens
c_p
 +  {\kappa p}w^{\mp 1}\tens (c_{p\mp 1}-c_p),
$$
where $a_\pm = a_1\pm \imath a_2$, and extended to the whole of
$E_\kappa(2)$ by 
(\ref{ent.A}). One way of 
seeing that the map $\psi$ is well-defined is to use the standard
coproduct of $E_\kappa(2)$ \cite{BonCic:new} and observe that  for any
$u\in E_\kappa(2)$, $\psi(c_l\tens u) = u\sw 1\tens \pi(w^lu\sw 2)$,
where $\pi : E_\kappa(2)\to C$ is a right module surjection given by
$\pi(w^l) =c_l$, $\pi(a_1) =\pi(a_2) =0$. By
\cite[Example~2.5]{BrzMa:coa} such a map necessarily entwines
$E_\kappa(2)$ with $C$. 
This can also be verified directly.
The map $\psi^C$ is given by
$
\psi^C(c_p\tens c_r) = c_r\tens c_{p+r}.
$

Using the map $\psi$ we define a right $C$-comodule structure on
$E_\kappa(2)$, $\Delta_R(u) = \psi(c_0\tens u)$. The fixed point
subalgebra $M$ of $E_\kappa(2)$ under this coaction is generated by
$a_1$, $a_2$ and  therefore it is isomorphic to the
$\kappa$-plane ${\bf R}_\kappa^2$. 

Finally, the crossed product structure of $E_\kappa(2)$ is provided by
the map 
$\Phi : C\to E_\kappa(2)$ given by
$\Phi(c_p) = v^{p}$. This map is clearly convolution invertible,
$\Phi(c_0) =1$ and the condition (\ref{triv}) is satisfied. Therefore
we have the crossed product $E_\kappa(2)= {\bf R}_\kappa^2\cross_{\rho,
\sigma} C$ 
 with a trivial cocycle $\sigma (c_p\tens c_r) =1$ and the map $\rho
:C\tens E_\kappa(2)\to E_\kappa(2)$, given explicitly by   
$$
\rho(c_p\tens  w^{\pm 1}) =1, \quad \rho(c_p\tens a_i) =a_i.
$$
Also in this case it is impossible to interpret $E_\kappa(2)$ as a
$C$-comodule (braided) algebra, and thus $E_\kappa(2)$ is truly an
example of a crossed product by a coalgebra.

\subsection{The $\kappa$-Euclidean group $E_\kappa(3)$}
Finally we sketch the crossed product structure of the deformation of
three-dimensional Euclidean group. $E_\kappa(3)$ is obtained by a
contraction of 
$SO_q(4)$, and  is generated by the Euler angles $\alpha , \beta, \gamma$, and
`coordinates' $x_1,x_2,x_3$ which satisfy the following relations
\cite{CelGia:...}: 
\begin{equation}
[\beta,x_+] = \kappa \sin\beta \tan(\beta/2), \quad [\beta,
x_3] = \kappa \sin\beta, \quad [x_-,\omega ] =
2\kappa\tan(\beta/2), 
\label{e31}
\end{equation}
\begin{equation}
[x_3,x_+] =\kappa x_3, \quad [x_-,x_+] = \kappa
{x_-},
\label{e32}
\end{equation}
where $x_+ =x_1\cos\alpha +x_2\sin\alpha$, $ x_- =
-x_1\sin\alpha +x_2\cos\alpha$, $\omega = \alpha +\gamma$ and all other
commutators vanish.  
The generators
$x_\pm,z$ span subalgebra of $E_\kappa(3)$ which is
denoted by $\R_\kappa^3$.
$E_\kappa(3)$ is a Hopf algebra and
$\R_\kappa^3$ is its homogeneous space, 
thus we  use \cite[Example~2.5]{BrzMa:coa} to
construct the entwining structure on $E_\kappa(3)$. Define a
right ideal $J\subset E_\kappa(3)$ generated by
$x_\pm,z$ and consider $C = E_\kappa(3)/J$. $C$ is
a coalgebra which can be equipped with a Hopf algebra structure of
functions on $SO(3)$ and is generated by the rotation matrix $A$ expressed
in terms of the Euler angles $\alpha,\beta,\gamma$. The
entwining structure is defined  by $\psi(c\tens u) = u\sw
1\otimes \pi(j(c)u\sw 2)$, where $\pi$ is a natural surjection
$E_\kappa(3) \to C$ and $j$ is a (Hopf algebra) inclusion
$C\hookrightarrow 
E_\kappa(3)$. Denoting elements of $A$
by $a_{ij}$, and using the fact that both $E_\kappa(3)$ and $C$ are
matrix quantum groups, 
we thus find 
\begin{equation}
\psi(a_{ij}\tens x_k) = x_k\tens a_{ij} + \sum_n a_{kn}\tens [a_{ij},x_n],
\quad \psi(a_{ij}\tens a_{kl}) = \sum_n a_{kn}\tens a_{ij}a_{nl},
\label{psi.e3}
\end{equation}
where the $SO(3)$ subgroup of $E_\kappa(3)$ is identified with $C$ via
$\pi$ and $j$. Equation (\ref{psi.e3}) can be made more explicit by
representing 
$a_{ij}$ in terms of the Euler angles and using (\ref{e31}-\ref{e32}).
% More explicitly, 
% \begin{eqnarray*}
% \psi(f(\alpha)g(\beta)\tens x_3) \!\!\! &=&\!\!\! x_3 \tens
% f(\alpha)g(\beta) + \kappa \cos\beta\tens
% f(\alpha)g'(\beta)\sin\beta\\ &&-\kappa (\sin\beta\cos\gamma \tens
% (f'(\alpha)g(\beta)\sin\alpha +
% f(\alpha)g'(\beta)\cos\alpha\sin\beta)\\
% &&+\sin\beta\sin\gamma\tens (f'(\alpha)g(\beta) \cos\alpha -
% f(\alpha)g'(\beta)\sin\alpha\sin\beta))\tan({\beta\over 2}) 
% \end{eqnarray*}
Furthermore we take 
$\psi^C(b\tens c) = c\sw 1\otimes bc\sw 2$ and $e=1$. Since $\pi$ is not an
algebra map, the coaction $\Delta_R(u) = 
\psi(1\tens u)$ is not an algebra map either, and to interpret
$E_\kappa(3)$ as a crossed product built on $\R_\kappa^3\tens C$, one
truly needs the construction 
described in Section~2. This is  provided by  the map $\Phi \equiv j$,
which clearly 
satisfies all the required conditions and thus makes $E_\kappa(3)$
into a crossed product $\R_\kappa^3\cross_{\rho,\sigma} C$ with a
trivial cocycle $\sigma(b,c) =\eps(b)\eps(c)$ and the 'action'
$$
\rho(c, x_i) = \eps(c)x_i, \qquad \rho(b,c) = \eps(b)\eps(c),
$$
for any $b,c\in C$.

\section{Conclusions}
In this short note we described crossed product structure of three
examples of inhomogeneous quantum groups. It is clear that other
inhomogeneous 
quantum groups may enjoy the similar structure, too (for example,
formula (\ref{psi.e3}) suggests an obvious generalisation).  We think that 
analysis of inhomogeneous quantum groups from the point of view of their
crossed product structure may enhance our understanding of these
groups, and thus indicate new applications to physics, special
functions and algebra. It can
also lead to development of gauge theories on such quantum groups
along the lines of \cite{BrzMa:coa}.

\end{document}